\documentclass[11pt,twoside]{article}
\usepackage{graphicx}
\usepackage{jeffe,smallbib,euler}
\usepackage{myconcrete,cyrillic}
\usepackage[same,tilde]{url}

\newtheorem{lemma}{Lemma}
\newtheorem{theorem}[lemma]{Theorem}

\begin{document}

\urldef{\paperURL}\url{http://www.cs.uiuc.edu/~jeffe/pubs/crum.html}

\title{\bf Arbitrarily Large Neighborly Families of\\
	Congruent Symmetric Convex 3-Polytopes}
\author{Jeff Erickson\thanks{Partially supported by a Sloan Research
	Fellowship and NSF CAREER award CCR-0093348.  See
	\paperURL\ for the most recent version of this paper.}
	\\[0.5ex]
	\small\begin{tabular}{c}
		University of Illinois at Urbana-Champaign\\
		\url{jeffe@cs.uiuc.edu}\\
		\url{http://www.cs.uiuc.edu/~jeffe}
	      \end{tabular}}
\maketitle
\thispagestyle{empty}

\begin{abstract}
We construct, for any positive integer $n$, a family of $n$ congruent 
convex polyhedra in $\Real^3$, such that every pair intersects in a 
common facet.  Previously, the largest such family contained only 
eight polytopes.  Our polyhedra are Voronoi regions of evenly 
distributed points on the helix $(t, \cos t, \sin t)$.  With a simple 
modification, we can ensure that each polyhedron in the family has a 
point, a line, and a plane of symmetry.  We also generalize our 
construction to higher dimensions and introduce a new family of cyclic 
polytopes.
\end{abstract}


\section{Introduction and History}

A family of $d$-dimensional convex polytopes is \emph{neighborly} if
every pair of polytopes has a $(d-1)$-dimensional intersection.  It
has been known for centuries that a neighborly family of convex
polygons (or any other connected sets) in the plane has at most four
members.  In 1905, Tietze \cite{t-udpdn-05,t-fpmsu-65} proved that
there are arbitrarily large neighborly families of $3$-dimensional
polytopes, answering an open question of Guthrie \cite{g--80} and
\Stackel~\cite{s--97}.  Tietze's result was independently rediscovered
by Besicovitch \cite{b-cp-47}, using a different construction, and
generalized to higher dimensions by Rado~\cite{e-recp-53} and
Eggleston \cite{r-sphis-47}.

Neighborly families of convex bodies are closely related to neighborly
convex polytopes.  A~polytope is ($2$-)neighborly if every pair of
vertices lies on a convex hull edge; the Schlegel diagram of the polar
dual of any neighborly $4$-polytope consists of a neighborly family of
$3$-polytopes.  Neighborly polytopes were discovered by
\Caratheodory~\cite{c-ud5df-11}, who showed that the convex hull of
any finite set of points on either the moment curve $(t, t^2, t^3,
\dots, t^d) \in \Real^d$ or the trigonometric moment curve $(\cos t,
\sin t, \allowbreak \cos 2t, \sin 2t, \dots, \allowbreak \cos kt, \sin
kt) \in \R^{2k}$ is a neighborly polytope.  \Caratheodory's proof was
simplified by Gale \cite{g-ncp-63}, who called these polytope families
the \emph{cyclic polytopes} and the \emph{Petrie polytopes},
respectively, and showed that the two families are combinatorially
equivalent.  Cyclic polytopes were independently rediscovered by
Motzkin \cite{m-ccp-57,gm-pg-63} and \Saskin~\cite{s-ravcp-63+}, among
others.  For further discussion of neighborly and cyclic polytopes,
see \Grunbaum~\cite{g-cp-67} and Ziegler~\cite{z-lp-94}.

Dewdney and Vranch \cite{dv-cprac-77} showed that the Voronoi diagram
of the integer points $\set{(t, t^2, t^3) \mid\allowbreak t =
1,2,\dots,n}$ form a neighborly family of unbounded convex polyhedra.
Klee \cite{k-cddvd-80} derived a similar result for any set of evenly
distributed points on the trigonometric moment curve in even
dimensions $4$ and higher.  Seidel \cite{s-eubnf-91} observed that for
any $d\ge 3$, Descartes' rule of signs\footnote{The number of real
roots of a polynomial is no more than the number of sign changes in
its degree-ordered sequence of non-zero coefficients.} implies that
any finite set of points on the positive branch of the $d$-dimensional
(polynomial) moment curve has a neighborly Voronoi diagram.  More
generally, the vertices of any neighborly polytope have a neighborly
Voronoi diagram, since the endpoints of any polytope edge have
neighboring Voronoi regions.

Zaks \cite{z-alnfs-86} described a general procedure to modify any 
neighborly family of unbounded polyhedra of any dimension, where each 
polyhedron contains an unbounded circular cone, so that the resulting 
polytopes are symmetric about a flat of any prescribed dimension.

Danzer, \Grunbaum, and Klee \cite{dgk-htir-63} asked if there is a 
largest neighborly family of \emph{congruent} polytopes.  Zaks (with 
Linhart) \cite{z-alnfs-86} observed that Klee's Voronoi diagram of 
evenly distributed points on the trigonometric moment curve forms a 
neighborly family of congruent convex polyhedra in even dimensions 
four and higher, but left the three-dimensional case open.  
The largest previously published neighborly family of congruent 
$3$-polytopes, discovered by Zaks~\cite{z-nfccp-87}, consists of 
eight triangular prisms.  According to Croft, Falconer, and Guy 
\cite[Problem E7]{cfg-upg-90}, this was also the largest known 
collection of congruent $3$-polytopes with the property that any two 
have even one point of contact.  Both Zaks~\cite{z-nfccp-87} and 
Croft, Falconer, and Guy \cite{cfg-upg-90} conjectured that the 
largest neighborly family of congruent $3$-polytopes is finite (see 
also Moser and Pach \cite[Problem 55]{mp-rpdgp-93}).

In Section 2, we show that this conjecture is incorrect, by giving a
constructive proof of the following theorem.

\begin{oneshot}{Main Theorem}
For any positive integer $n$, there is a neighborly family of $n$
congruent convex $3$-polytopes.
\end{oneshot}

\begin{figure}[htb]
\centerline{\footnotesize\sf
\begin{tabular}{c@{\qquad}c}
   \includegraphics[height=2in]{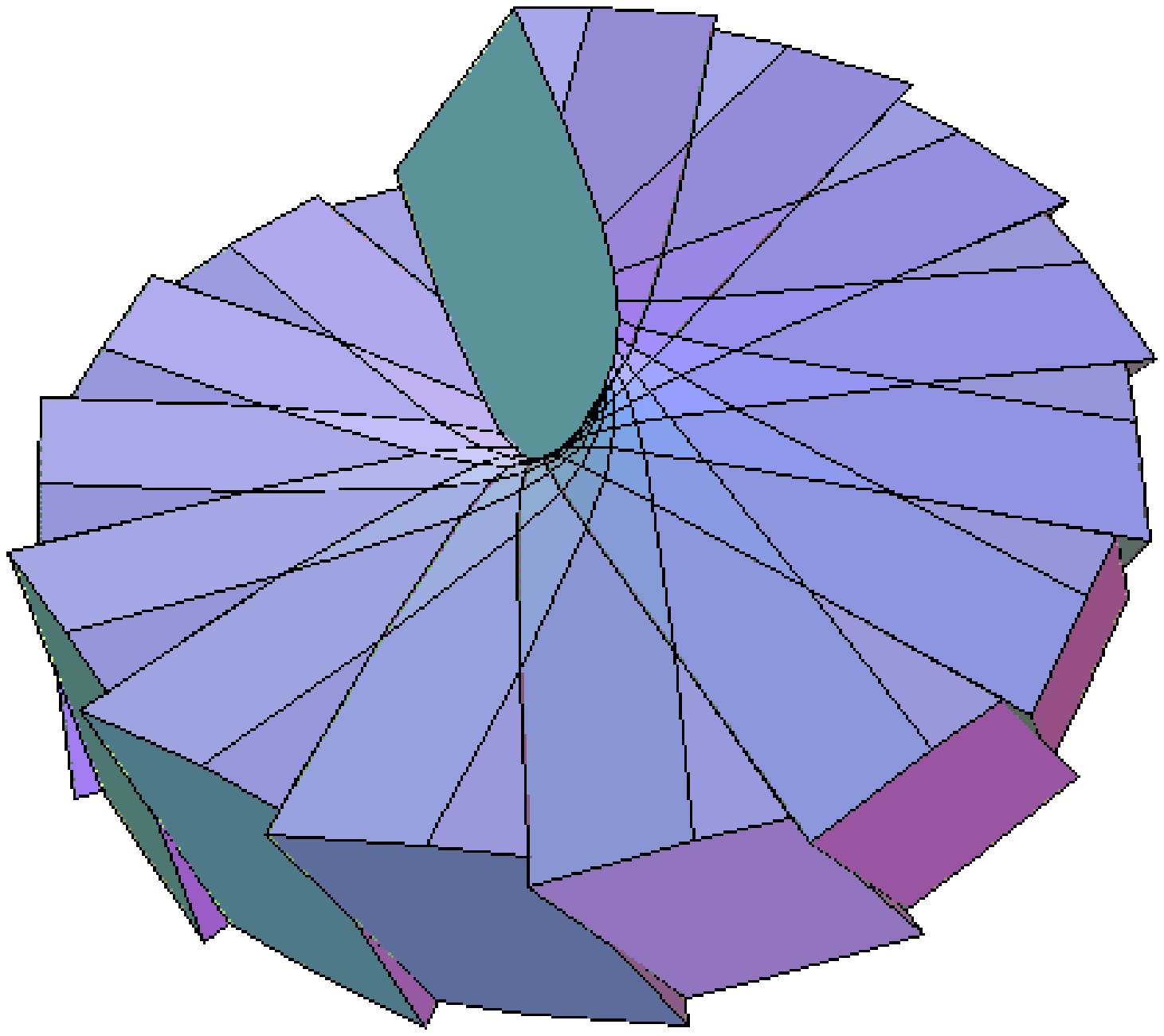} &
   \includegraphics[height=2in]{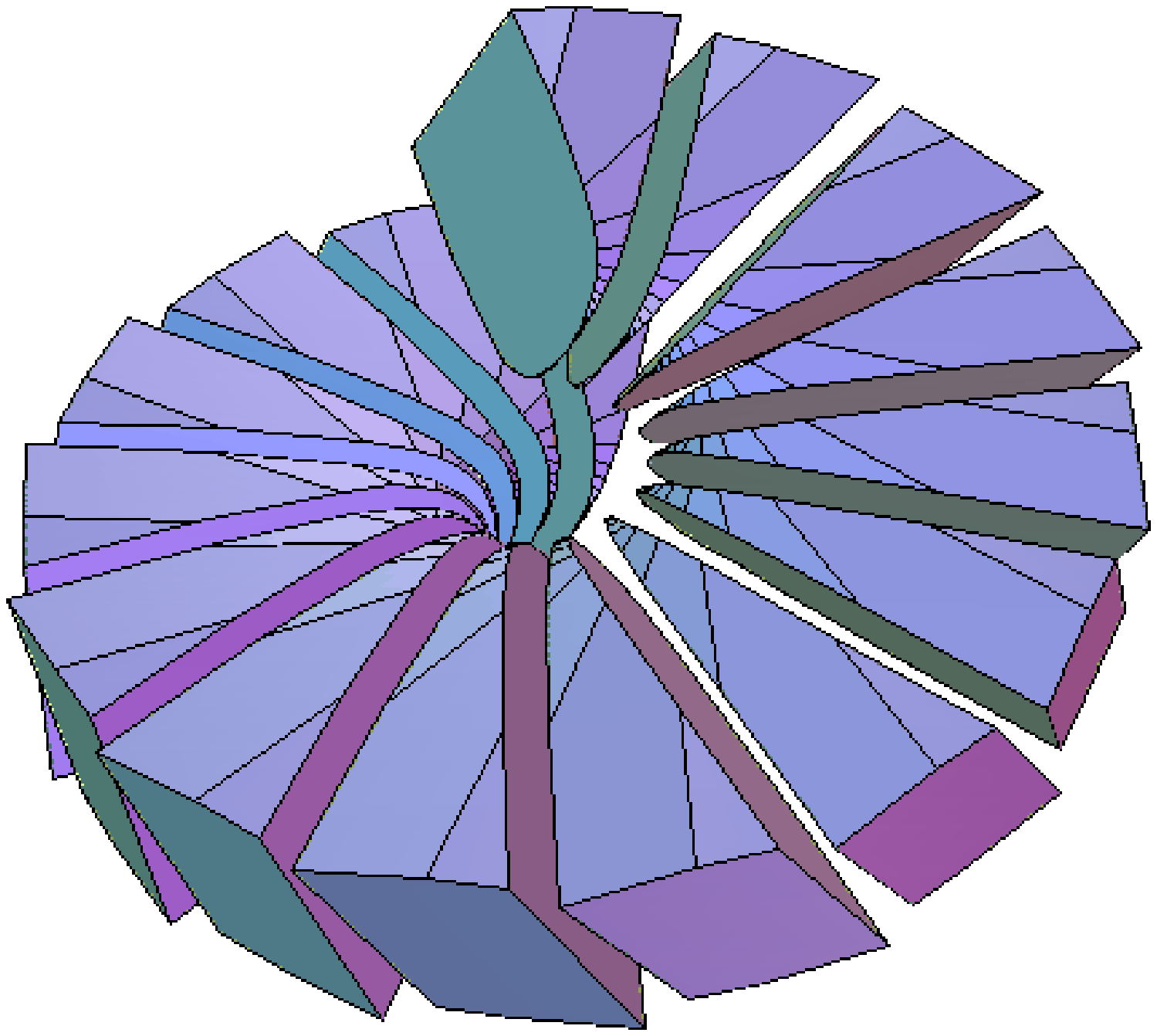}\\
(a) & (b)
\end{tabular}}
\caption{(a)~A neighborly family of sixteen congruent convex
polytopes.  (b)~An exploded view of the same family.}
\label{Fig/snail}
\end{figure}

Like the earlier constructions of Dewdney and 
Vranch~\cite{dv-cprac-77} and Zaks and Linhart \cite{z-nfccp-87}, our 
construction is based on the Voronoi diagram of a set of points on a 
curve, namely the regular circular helix $h(t) = (t, \cos t, \sin t)$.  
An example of our construction is shown in Figure~\ref{Fig/snail}, and 
a single polytope in our family is shown in Figure \ref{Fig/snail1}.

We generalize our Main Theorem to higher dimensions in 
Section~\ref{S:higher}, by constructing an arbitrarily large family of 
congruent convex polytopes in $\Real^d$, any $\ceil{d/2}$ of which 
share a unique common boundary face.  We also introduce a new family 
of cyclic polytopes, generalizing both the classic cyclic polytopes 
and the Petrie polytopes.

\section{The Main Theorem}
\label{S:main}

Our construction relies on the following observation of the author
\cite{e-npsch-01}.  We include the proof for the sake of completeness.

\begin{lemma}
\label{L:bitangent}
Let $\beta(t)$ denote the unique sphere passing through $h(t)$ and 
$h(-t)$ and tangent to the helix at those two points.  For any 
$0<t<pi$, the sphere $\beta(t)$ intersects the helix only at its two 
points of tangency.
\end{lemma}

\begin{proof}
Since a $180$-degree rotation about the $y$-axis maps $h(t)$ to $h(-t)$
and leaves the helix invariant, the bitangent sphere must be centered
on the $y$-axis.  Thus, $\beta$ can be described by the equation $x^2
+ (y-a)^2 + z^2 = r^2$ for some constants $a$ and~$r$.  Let $\gamma$
denote the intersection curve of $\beta(t)$ and the cylinder $y^2 +
z^2 = 1$.  Every intersection point between $\beta(t)$ and the helix
must lie on $\gamma$.  If we project the helix and the intersection
curve to the $xy$-plane, we obtain the sinusoid $y = \cos x$ and a
portion of the parabola $y = \gamma(x) = (x^2 - r^2 + a^2 + 1)/2a$.
These two curves meet tangentially at the points $(t, \cos t)$ and
$(-t, \cos t)$.

\begin{figure}[htb]
    \centerline{\includegraphics[height=1in]{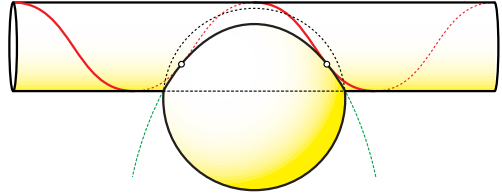}}
    \caption{The intersection curve of the cylinder and a bitangent 
    sphere projects to a parabola on the $xy$-plane.}
\end{figure}

The mean value theorem implies that $\gamma(x) = \cos x$ at most four
times in the range $-\pi < x < \pi$.  (Otherwise, the curves $y'' =
-\cos x$ and $y'' = \gamma''(x) = 1/a$ would intersect more than twice
in that range.)  Since the curves meet with even multiplicity at two
points, those are the only intersection points in the range $-\pi < x
< \pi$.  Since $\gamma(x)$ is concave, we have $\gamma(\pm\pi) <
\cos\pm\pi = -1$, so there are no intersections with $\abs{x}\ge\pi$.
Thus, the curves meet only at their two points of tangency.
\end{proof}

Lemma \ref{L:bitangent} immediately implies that the Voronoi diagram
of any finite set of points on the helix $(t, \cos t, \sin t)$ in the
range $-\pi < t < \pi$ is a neighborly collection of unbounded convex
polyhedra.

\def\Hn{\mathscr{H}_n}

To obtain a neighborly family of \emph{congruent} polyhedra, we use 
the Voronoi diagram of evenly spaced points on the helix.  For any 
integer $n$, let $h_n(t) = h(2\pi t/n)$ and let $\Hn$ denote the 
infinite point set $\set{h_n(t) \mid t\in\Z}$.  By the Bitangent 
Lemma, the Voronoi regions of any $n+1$ consecutive points in $\Hn$ 
form a neighborly family of convex bodies.  Since the point set 
$\mathscr{H}_n$ is preserved by the rigid motion
\[
	(x,y,z) \mapsto 
	\left(	x+\frac{2\pi}{n},\;
		y\,\cos\frac{2\pi}{n} - z\,\sin\frac{2\pi}{n},\;
		y\,\sin\frac{2\pi}{n} + z\,\cos\frac{2\pi}{n}
	\right),
\]
which maps each point $h_n(t)$ to its successor $h_n(t+1)$, these 
Voronoi regions are all congruent.

The following more refined analysis of the Delaunay triangulation of 
$\Hn$, reminiscent of Gale's `evenness condition' for cyclic 
polytopes~\cite{g-ncp-63,s-tcp-68}, implies that these Voronoi regions 
have only a finite number of facets, and thus are actually polyhedra.

\begin{lemma}\label{L:pairs}
For any integers $a < b < c < d$, the points $h_n(a), h_n(b), h_n(c),
h_n(d)$ are vertices of a simplex in the Delaunay triangulation of
$\Hn$ if and only if $b - a = d - c = 1$ and $d - a \le n$.
\end{lemma}

\begin{proof}
Call a tetrahedron with vertices $h_n(a), h_n(b), h_n(c), h_n(d)$
\emph{local} if $b - a = d - c = 1$ and $d - a \le n$.  Let $\sigma$
be the sphere passing through the vertices of an arbitrary local
tetrahedron.  Analysis similar to the proof of Lemma~\ref{L:bitangent}
implies that the only portions of the helix that lie inside $\sigma$
are the segments between $h_n(a)$ and $h_n(b)$ and between $h_n(c)$
and $h_n(d)$.  Thus, all other points in $\Hn$ lie outside $\sigma$,
so the four points form a Delaunay simplex.

The local Delaunay simplices exactly fill the convex hull of $\Hn$,
and therefore comprise the entire Delaunay triangulation.
Specifically, the only triangles that are facets of exactly one local
tetrahedron have vertices $h_n(i), h_n(i+1), h_n(i+n)$ or $h_n(i-n+1),
h_n(i), h_n(i+1)$ for some integer $i$.  Thus, a tetrahedron is
Delaunay if and only if it is local.
\end{proof}

\begin{figure}[htb]
\centerline{\footnotesize\sf
\begin{tabular}{c@{\qquad}c}
   \includegraphics[width=1.25in,angle=-90]{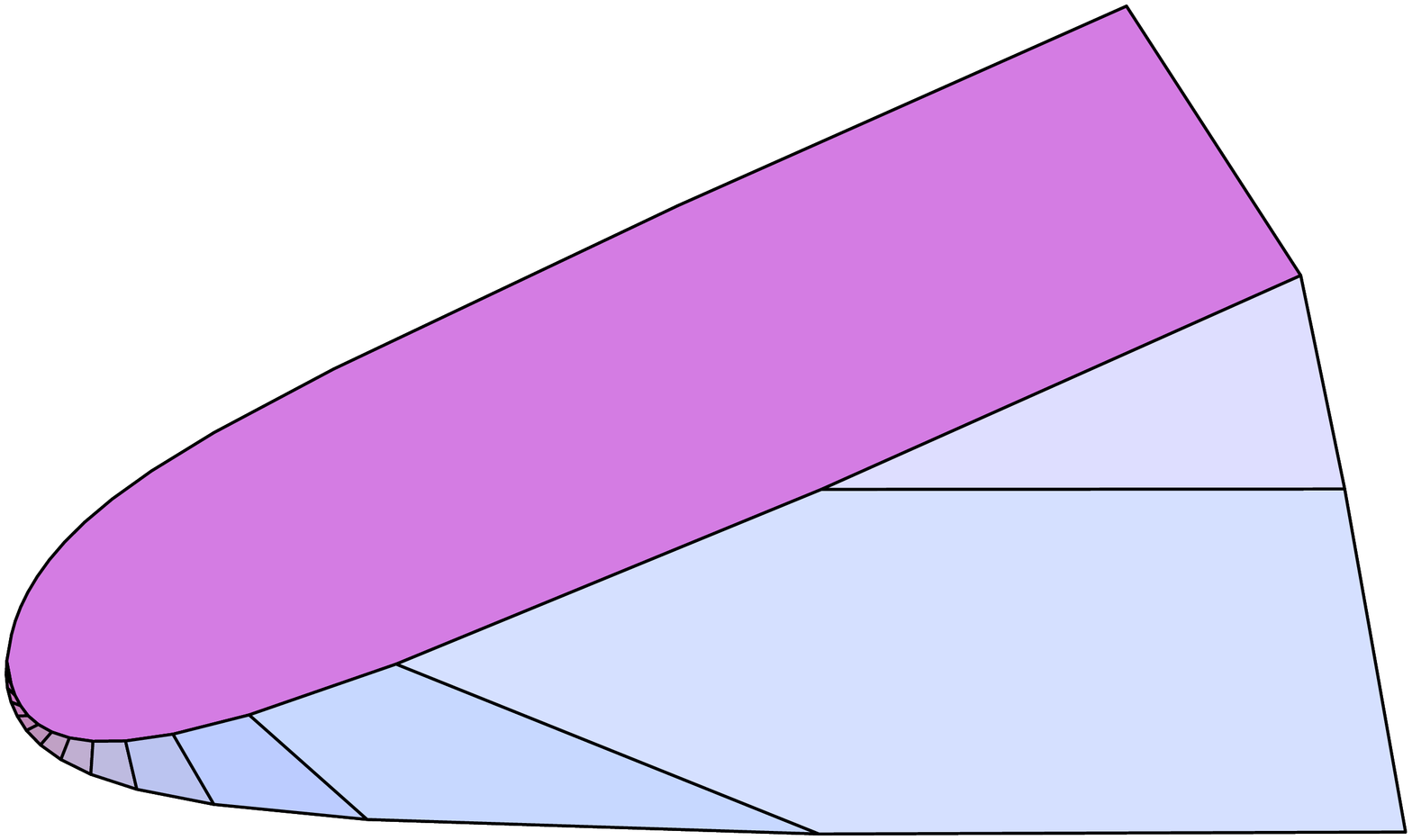} & 
   \includegraphics[width=1.25in,angle=-90]{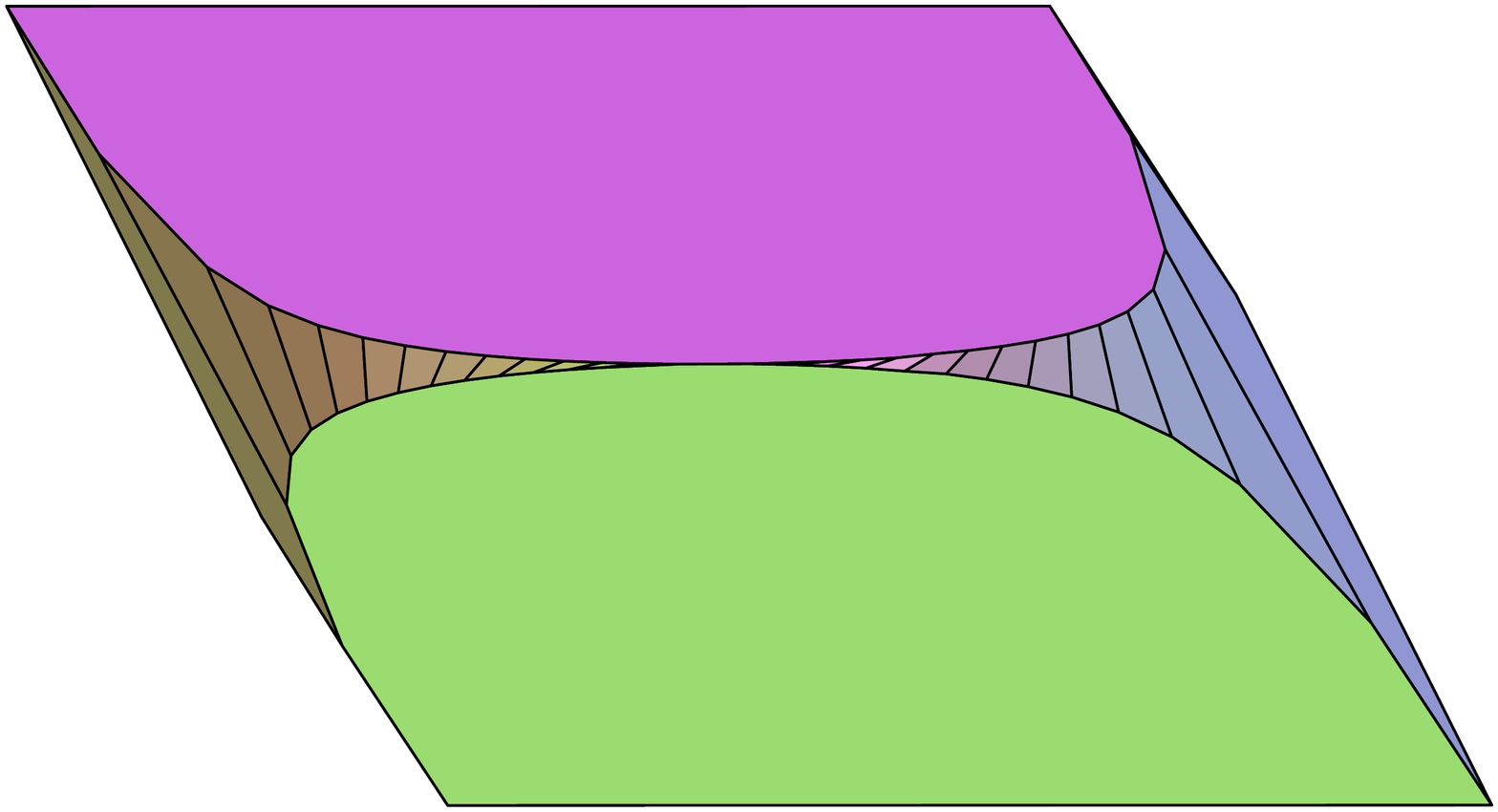} \\
 (a) & (b)
\end{tabular}}
\caption{(a) One polytope in the neighborly family of sixteen.  (b)~An
orthographic edge-on view of the same polytope.}
\label{Fig/snail1}
\end{figure}

In light of the duality between features of Delaunay triangulations
and Voronoi diagrams, Lemma~\ref{L:pairs} lets us exactly describe the
combinatorial structure of the Voronoi regions of~$\Hn$.  Let $V_n(t)$
denote the Voronoi region of $h_n(t)$.  This polyhedron has exactly
$2n$ facets, in $n$ symmetric pairs, as follows:
%
%
\begin{itemize}\itemsep=0pt
\item[$\bullet$]
two unbounded $(2n-1)$-gons shared with $V_n(t\pm 1)$, each bounded by
$2n-3$ segments and two parallel rays;
\item[$\bullet$]
two triangles shared with $V_n(t\pm 2)$;
\item[$\bullet$]
$2n-8$ quadrilaterals shared with $V_n(t\pm 3), V_n(t\pm 4), \dots,
V_n(t\pm (n-2))$;
\item[$\bullet$]
two unbounded quadrilaterals shared with $V_n(t\pm (n-1))$,
each bounded by two line segments and two parallel rays;
\item[$\bullet$]
two wedges in parallel planes shared with $V_n(t\pm n)$, each bounded
by a pair of rays.
\end{itemize}
The two $(2n-1)$-gons are adjacent to all the other facets, including 
each other, and contain all the vertices of $V_n(t)$; otherwise, the 
facets are adjacent in sequence.  See Figure \ref{Fig/snail1}.

The vertex of $V_n(0)$ furthest from the $x$-axis is the center of the 
sphere through $h_n(0)$, $h_n(1)$, $h_n(n-1)$ and $h_n(n)$, which has 
coordinates $(\pi, \theta(2\pi-\theta)/(2-2\cos\theta), 0)$, where 
$\theta = 2\pi/n$.  Thus, all the Voronoi vertices of $\Hn$ lie in a 
cylinder of radius $\theta(2\pi-\theta)/(2-2\cos\theta) \approx n-1$ 
around the $x$-axis.  To transform our neighborly family of unbounded 
polyhedra into a neighborly family of \emph{polytopes}, we intersect 
each Voronoi region $V_n(t)$ with the halfspace $y \cos (2\pi t/n) + z 
\sin (2\pi t/n) \le n$, which contains some positive area of every 
facet of $V_n(t)$.

This completes the proof of the Main Theorem.

\medskip Zaks \cite{z-pc-01} describes an alternate proof, based
entirely on the neighborliness of the Voronoi regions of~$\Hn$.  For
each integer $1\le t\le n$, place a triangle on the shared boundary
facet between $V_n(0)$ and $V_n(t)$.  Now place congruent copies of
these triangles on the boundary of every Voronoi region, so that the
entire collection has the same screw symmetry as~$\Hn$.  Finally, for
any integer~$t$, let $C_n(t)$ be the convex hull of the triangles on
the boundary of $V_n(t)$.  The $n+1$ congruent convex polytopes
$C_n(0), C_n(1), \dots, C_n(n-1), C_n(n)$ form a neighborly family.

To actually construct either our neighborly family or Zaks', it
suffices to compute the Voronoi diagram of the finite point set
$\set{h_n(t) \mid t=0,1,2,\dots,3n}$ and then consider only the
Voronoi regions of the middle $n+1$ points $h_n(n), h_n(n+1), \dots,
h_n(2n-1), h_n(2n)$, since those Voronoi regions are the same as in
the infinite point set $\Hn$.  Figure 1 was computed using this
method.

Finally, since a $180$-degree rotation about the $y$-axis maps each 
point $h_n(t)$ to $h_n(-t)$, and thus preserves the point set $\Hn$, 
the Voronoi region $V_n(0)$ is rotationally symmetric about the 
$y$-axis.  It immediately follows every Voronoi region of $\Hn$ has a 
line of $180$-degree rotational symmetry.  Clipping each Voronoi 
region by an additional halfspace as above retains this symmetry, sine 
the clipping plane is normal to the symmetry axis.  We can create a 
neighborly family of congruent polytopes additional symmetries by 
taking the union of each clipped Voronoi region and its reflection 
across its clipping plane.  Each resulting polytope clearly has 
bilateral symmetry about its clipping plane and $180$-degree symmetry 
about the original Voronoi region's axis of symmetry, and therefore is 
centrally symmetric about the intersection point of the clipping plane 
and the symmetry axis.

\begin{theorem}
For any integer positive integer $n$, there is a neighborly family of 
$n$ congruent convex $3$-polytopes, each with a plane of bilateral 
symmetry, a line of $180$-degree rotational symmetry, and a point of 
central symmetry.
\end{theorem}

\section{Higher Dimensions}
\label{S:higher}

A family of convex polyhedra in $\Real^d$ is \emph{(strictly)
$k$-neighborly} if any subset of $k$ polyhedra has a
$(d-k+1)$-dimensional intersection, and no subset of $k+1$ polytopes
has a non-empty intersection.\footnote{The second condition is
necessary to rule out degenerate constructions such as the product of
a $(d-2)$-dimensional cube with $n$ congruent planar wedges.}
Arbitrarily large $k$-neighborly families of polyhedra are easy to
construct in $\Real^{2k-1}$, for example, Schlegel diagrams of dual
cyclic $2k$-polytopes \cite{c-ud5df-11,g-ncp-63} or Voronoi diagrams
of points on the moment curve \cite{s-eubnf-91}.  However, arbitrarily
large $k$-neighborly families of \emph{congruent} polyhedra were
previously only known in dimensions $2k$ and higher.  The
lowest-dimensional example is based on the Voronoi diagram of evenly
distributed points on the trigonometric moment curve
\cite{k-cddvd-80,z-alnfs-86} together with the origin (since otherwise
the origin is on the boundary of every Voronoi polyhedron).

In this section, we generalize our three dimensional results by 
considering regularly spaced points on the following \emph{generalized 
helix}:
\[
	h_k(t) =
	\big(	t,\;
		\cos t,\, \sin t,\,
		\cos 2t,\, \sin 2t,\, \dots,\,
		\cos kt,\, \sin kt
		\big)
	\in \Real^{2k+1}.
\]

\begin{theorem}
\label{Th:oneturn}
Let $P$ be any finite set of points on the curve $h_k(t)$ in the range
$0 < t < 2\pi$, for some non-negative integer $k$.  The Voronoi
diagram of $P$ is a $(k+1)$-neighborly family of convex polyhedra in
$\Real^{2k+1}$.
\end{theorem}

\begin{proof}
Consider the sphere $\sigma$ passing through $k+1$ arbitrary points
$h_k(a_0), h_k(a_1), \dots,\allowbreak {h_k(a_k) \in P}$ and tangent
to the generalized helix at those points, where $0 < a_0 < a_1 <
\cdots < a_k < 2\pi$.  Any point $h_k(t)$ that lies on $\sigma$
satisfies the following matrix equation:
\[
F(t) = 
\begin{vmatrix}
	1 & a_0 & \cos a_0   & \sin a_0
		& \cos 2 a_0 & \sin 2 a_0
		& \cdots
		& \cos k a_0 & \sin k a_0
		& k + a_0^2
\\
	0 & 1   & -\sin a_0     & \cos a_0
		& -2\sin 2 a_0  & 2\cos 2 a_0
		& \cdots
		& -k \sin k a_0 & k\cos k a_0
		& 2 a_0
\\
	1 & a_1 & \cos a_1   & \sin a_1
		& \cos 2 a_1 & \sin 2 a_1
		& \cdots
		& \cos k a_1 & \sin k a_1
		& k + a_1^2
\\
	0 & 1   & -\sin a_1     & \cos a_1
		& -2\sin 2 a_1  & 2\cos 2 a_1
		& \cdots
		& -k \sin k a_1 & k\cos k a_1
		& 2 a_1
\\
	\vdots  & \vdots
		& \vdots & \vdots
		& \vdots & \vdots
		& \ddots
		& \vdots & \vdots
		& \vdots
\\
	1 & a_k & \cos a_k   & \sin a_k
		& \cos 2 a_k & \sin 2 a_k
		& \cdots
		& \cos k a_k & \sin k a_k
		& k + a_k^2
\\
	0 & 1   & -\sin a_k     & \cos a_k
		& -2\sin 2 a_k  & 2\cos 2 a_k
		& \cdots
		& -k \sin k a_k & k\cos k a_k
		& 2 a_k
\\
	1 & t   & \cos t     & \sin t
		& \cos 2 t   & \sin 2 t
		& \cdots
		& \cos k t   & \sin k t
		& k + t^2
\end{vmatrix} = 0
\]
To bound the number of zeros of $F(t)$, consider its second derivative
\[
F''(t) =
\begin{vmatrix}
	1 & a_0 & \cos a_0   & \sin a_0
		& \cos 2 a_0 & \sin 2 a_0
		& \cdots
		& \cos k a_0 & \sin k a_0
		& k + a_0^2
\\
	0 & 1   & -\sin a_0     & \cos a_0
		& -2\sin 2 a_0  & 2\cos 2 a_0
		& \cdots
		& -k \sin k a_0 & k\cos k a_0
		& 2 a_0
\\
	1 & a_1 & \cos a_1   & \sin a_1
		& \cos 2 a_1 & \sin 2 a_1
		& \cdots
		& \cos k a_1 & \sin k a_1
		& k + a_1^2
\\
	0 & 1   & -\sin a_1     & \cos a_1
		& -2\sin 2 a_1  & 2\cos 2 a_1
		& \cdots
		& -k \sin k a_1 & k\cos k a_1
		& 2 a_1
\\
	\vdots  & \vdots
		& \vdots & \vdots
		& \vdots & \vdots
		& \ddots
		& \vdots & \vdots
		& \vdots
\\
	1 & a_k & \cos a_k   & \sin a_k
		& \cos 2 a_k & \sin 2 a_k
		& \cdots
		& \cos k a_k & \sin k a_k
		& k + a_k^2
\\
	0 & 1   & -\sin a_k     & \cos a_k
		& -2\sin 2 a_k  & 2\cos 2 a_k
		& \cdots
		& -k \sin k a_k & k\cos k a_k
		& 2 a_k
\\
	0 & 0   & -\cos t     & -\sin t
		& -4\cos 2 t  & -4\sin 2 t
		& \cdots
		& -k^2\cos kt & -k^2 \sin k t
		& 2
\end{vmatrix}.
\]
$F''(t)$ is an affine combination of the functions $\cos t, \sin t,
\allowbreak \cos 2t, \sin 2t, \dots,\allowbreak \cos kt, \sin kt$, so
it can be rewritten as as a polynomial of degree at most $2k$ in the
variable $e^{it}$.  Thus, $F''(t)$ has at most $2k$ zeros in the range
$0 < t < 2\pi$.  (This is essentially the argument used by
\Caratheodory\ to show that Petrie polytopes are neighborly
\cite{c-ud5df-11}.)

Since $a_0, a_1, \dots, a_k$ are roots of $F(t)$ of multiplicity two,
they are the only roots in the range $0 < t < 2\pi$; otherwise, by the
mean value theorem, $F''(k)$ would have more than $2k$ roots in the
range $0\le a_0 < t < a_k\le 2\pi$, which we have just shown to be
impossible.  Thus, the points $h_k(a_0),\allowbreak h_k(a_1),
\dots,\allowbreak h_k(a_k)$ lie on a sphere that excludes every other
point in $P$ and so have mutually neighboring Voronoi regions.
\end{proof}

In fact, Theorem \ref{Th:oneturn} is a special case of the following
result, which follows from an easy generalization of the previous
proof and Gale's evenness condition for cyclic polytopes
\cite{g-ncp-63,s-tcp-68}.  Define the \emph{mixed moment curve}
$\mu_{d,k}(t)$ as follows:
\[
	\mu_{d,k}(t) = \big(t, t^2, \dots, t^d,\;
			\cos t, \sin t,\,
			\cos 2t, \sin 2t,\, 
			\dots,
			\cos kt, \sin kt
			\big) \in \Real^{2k+d}.
\]
For example, $\mu_{d,0}(t)$ is the standard $d$-dimensional moment 
curve, $\mu_{0,k}(t)$ is the $2k$-dimensional trigonometric moment 
curve, and $\mu_{1,k}(t)$ is our generalized helix.

\begin{theorem}
For any non-negative integers $d$ and $k$, the convex hull of any
finite set of points on the curve $\mu_{d,k}(t)$ in the range $0 < t <
2\pi$ is a $(d+2k)$-dimensional cyclic polytope.
\end{theorem}

Not surprisingly, we obtain large highly-neighborly families of
\emph{congruent} polytopes by considering the Voronoi diagram of
infinite point set $\mathcal{H}_n^k = \set{h_k(2\pi t/n) \mid
t\in\Z}$.

\begin{theorem}
For any non-negative integers $n$ and $k$, any $n+1$ consecutive
Voronoi regions in the Voronoi diagram of $\mathcal{H}_n^k$ form a
$(k+1)$-neighborly family of congruent convex polyhedra.
\end{theorem}

\begin{proof}
Fix an integer $n$, and for notational convenience, let $\hbar(t) = 
h_k(2\pi t/n)$.  Since $\mathcal{H}^k_n$ is preserved under a rigid 
motion mapping each point $\hbar(i)$ to its successor $\hbar(i+1)$, the 
Voronoi regions of $\mathcal{H}^k_n$ are congruent.

Call a full-dimensional simplex with vertices in $\mathcal{H}^k_n$
\emph{local} if all its vertices consist of $k+1$ adjacent pairs
within a single turn of the generalized helix, that is, if its
vertices are
\[
	\hbar(a_0),\; \hbar(a_0+1),\; \hbar(a_1),\; 
	\hbar(a_1+1),\; \dots,\; \hbar(a_k),\; \hbar(a_k+1),
\]
for some integers $a_0, a_1, \dots, a_k$ with $a_k+1 \le a_0+n$ and
$a_i + 1 < a_{i+1}$ for all $i$.  Analysis similar to Theorem
\ref{Th:oneturn} implies that every local simplex is Delaunay.

The convex hull of $\mathcal{H}^k_n$, which we will call the
\emph{Petrie cylinder}, is the product of an $2k$-dimensional Petrie
polytope with $n$ vertices and a line orthogonal to that polytope's
hyperplane.  By Gale's evenness condition \cite{g-ncp-63,s-tcp-68},
the facets of the Petrie polytope are formed by all sets of $k$
adjacent pairs of points on the trigonometric moment curve.  The faces
of the Petrie cylinder are cylinders over the faces of the Petrie
polytope.

Call a facet of a local simplex that is not shared by another local
simplex a \emph{boundary simplex}.  We easily observe that the
boundary simplices are exactly the $2k$-simplices whose ordered
sequence of vertices has one of the following two forms:
\def\arraystretch{1.2}
\[
\begin{array}{r@{}l}
	\langle
	\hbar(a_k+1-n), \hbar(a_1), \hbar(a_1+1), \hbar(a_2),&
	\hbar(a_2+1), \dots, \hbar(a_k), \hbar(a_k+1)
	\rangle
\\
	\langle
	\hbar(a_1), \hbar(a_1+1), \hbar(a_2),&
	\hbar(a_2+1), \dots, \hbar(a_k), \hbar(a_k+1), \hbar(a_1+n)
	\rangle
\end{array}
\]
The following sequence of boundary simplices exactly covers one facet
of the Petrie cylinder.
{\small
\[
\begin{array}{r@{}l}
	\ddots\qquad\qquad\hphantom{\dots,}&
\\
	\langle
	\hbar(a_k+1-n), \hbar(a_1), \hbar(a_1+1), \hbar(a_2),
	\hbar(a_2+1), \dots,
	&
	\hbar(a_k), \hbar(a_k+1)
	\rangle,
\\
	\langle
	\hbar(a_1), \hbar(a_1+1), \hbar(a_2), \hbar(a_2+1), \dots,
	&
	\hbar(a_k), \hbar(a_k+1), \hbar(a_1+n)
	\rangle,
\\
	\langle
	\hbar(a_1+1), \hbar(a_2), \hbar(a_2+1), \dots,
	&
	\hbar(a_k), \hbar(a_k+1), \hbar(a_1+n), \hbar(a_1+n+1)
	\rangle,
\\
	\langle
	\hbar(a_2), \hbar(a_2+1), \dots,
	&
	\hbar(a_k), \hbar(a_k+1), \hbar(a_1+n), \hbar(a_1+n+1), \hbar(a_2+n)
	\rangle,
\\
	&\qquad\qquad\ddots
\end{array}
\]
}%
Every facet of the Petrie cylinder is covered in this manner, and
every boundary simplex lies on some facet of the Petrie cylinder.
Thus, the union of the boundary facets is the boundary of the Petrie
cylinder, so the local Delaunay simplices completely fill the Petrie
cylinder and therefore comprise the entire Delaunay triangulation.

It easily follows that each Voronoi region of $\mathcal{H}^k_n$ is a 
convex polyhedron with $\Theta(n^k)$ facets, and that any $n+1$ 
consecutive Voronoi regions form a $(k+1)$-neighborly family.  As we 
already observed, these polyhedra are congruent.
\end{proof}

We can easily modify our construction to obtain a $(k+1)$-neighborly
family of polytopes, by intersecting each Voronoi region with a
halfspace strictly containing all its vertices.  Each Voronoi region
of $\mathcal{H}^k_n$ has a $k$-flat of two-fold symmetry.  As long as
the boundary of the new halfspace is perpendicular to this central
$k$-flat, the resulting polytope is also symmetric about this flat.

Using a variant of Zaks' symmetrization procedure \cite{z-alnfs-86},
we can ensure that each polytope is also symmetric about a flat of any
specified dimension.  Consider the Voronoi region $V$ of $\hbar(0)$ in
the Voronoi diagram of $\mathcal{H}^k_n$.  Let $\rho$ be the ray from
the origin through $\hbar(0)$, let $\phi^+$ and $\phi^-$ denote the
supporting hyperplanes of the only two parallel facets of $V$ (shared
with the Voronoi regions of $\hbar(n)$ and $\hbar(-n)$), and let $\pi$
be a hyperplane normal to $\rho$ at sufficient distance from the
origin.  Finally, let $f$ be any flat that lies in $\pi$, contains the
point $\rho\cap\pi$, and is either parallel or perpendicular to
$\phi^+$ and~$\phi^-$.  The intersection of $V$ and its reflection
across $f$ is a convex polytope that is obviously symmetric about $f$,
and whose boundary contains positive measure from every boundary facet
of~$V$.  Applying this procedure to any $n+1$ consecutive Voronoi
regions of $\mathcal{H}^k_n$, we obtain our final result.

\begin{theorem}
For any positive integers $k$ and $n$ and any nonnegative integer
$r\le 2k$, there is a $(k+1)$-neighborly family of $n$ congruent
convex polytopes in $\Real^{2k+1}$, each of which is symmetric about
an $r$-flat.
\end{theorem}

\paragraph{Acknowledgments.} Thanks to Joseph Zaks for his insightful 
comments on an early draft of this paper, and to Victor Klee for 
sending me a copy of his paper \cite{k-cddvd-80}.  Figures 
\ref{Fig/snail} and \ref{Fig/snail1} were produced with the help of 
the programs 'qhull'~\cite{bdh-qach-96,qhull} and 
`geomview'~\cite{almp-gsgv-95,geomview}.

\bibliographystyle{abuser}
\bibliography{crum,geom}

\def\burl#1{$\langle$\url{#1}$\rangle$}
\begin{thebibliography}{10}

\bibitem{almp-gsgv-95}
N.~Amenta, S.~Levy, T.~Munzner, and M.~Philips.
\newblock Geomview: A system for geometric visualization.
\newblock \emph{Proc. 11th Annu. ACM Sympos. Comput. Geom.}, pp. C12--C13,
  1995.

\bibitem{bdh-qach-96}
C.~B. Barber, D.~P. Dobkin, and H.~Huhdanpaa.
\newblock The \textsc{Quickhull} algorithm for convex hulls.
\newblock \emph{ACM Trans. Math. Softw.} 22(4):469--483, Dec. 1996.

\bibitem{qhull}
C.~B. Barber and H.~Huhdanpaa.
\newblock Qhull, version 3.0, February 2000.
\newblock \burl{http://www.geom.umn.edu/software/qhull/}.

\bibitem{b-cp-47}
A.~S. Besicovitch.
\newblock On {C}rum's problem.
\newblock \emph{J. London Math. Soc.} 22:285--287, 1947.

\bibitem{c-ud5df-11}
C.~Carath{\'e}odory.
\newblock {\"U}ber den {V}ariabilit{\"a}tsbereich der {F}ourier'schen
  {K}onstanten von positiven harmonischen {F}unktionen.
\newblock \emph{Rendiconto del Circolo Matematico di Palermo} 32:193--217,
  1911.

\bibitem{cfg-upg-90}
H.~P. Croft, K.~J. Falconer, and R.~K. Guy.
\newblock \emph{Unsolved Problems in Geometry}.
\newblock Springer-Verlag, 1990.

\bibitem{dgk-htir-63}
L.~Danzer, B.~Gr{\"u}nbaum, and V.~Klee.
\newblock Helly's theorem and its relatives.
\newblock \emph{Convexity}, pp. 101--180. Proc. Symp. Pure Math. VII, Amer.
  Math. Soc., 1963.

\bibitem{dv-cprac-77}
A.~K. Dewdney and J.~K. Vranch.
\newblock A convex partition of {$R^{3}$} with applications to {Crum}'s problem
  and {Knuth}'s post-office problem.
\newblock \emph{Utilitas Math.} 12:193--199, 1977.

\bibitem{e-recp-53}
H.~G. Eggleston.
\newblock On {Rado's} extension of {Crum's} problem.
\newblock \emph{J. London Math. Soc.} 28:467--471, 1953.

\bibitem{e-npsch-01}
J.~Erickson.
\newblock Nice point sets can have nasty {D}elaunay triangulations.
\newblock \emph{Proc. 17th Annu. ACM Sympos. Comput. Geom.}, pp. 96--105,
  2001.
\newblock arXiv:cs.CG/0103017.

\bibitem{g-ncp-63}
D.~Gale.
\newblock Neighborly and cyclic polytopes.
\newblock \emph{Convexity}, pp. 225--232. Proc. Symp. Pure Math. VII, Amer.
  Math. Soc., 1963.

\bibitem{g-cp-67}
B.~Gr{\"u}nbaum.
\newblock \emph{Convex Polytopes}.
\newblock John Wiley \& Sons, New York, NY, 1967.

\bibitem{gm-pg-63}
B.~Gr{\"u}nbaum and T.~S. Motzkin.
\newblock On polyhedral graphs.
\newblock \emph{Convexity}, pp. 285--290. Proc. Symp. Pure Math. VII, Amer.
  Math. Soc., 1963.

\bibitem{g--80}
F.~Guthrie.
\newblock \emph{Proc. Royal Soc. Edinburgh} 10:728, 1878--1880.

\bibitem{k-cddvd-80}
V.~Klee.
\newblock On the complexity of $d$-dimensional {V}oronoi diagrams.
\newblock \emph{Archiv der Math.} 34:75--80, 1980.

\bibitem{geomview}
S.~Levy, T.~Munzner, M.~Phillips, C.~Fowler, N.~Thurston, D.~Krech, S.~Wisdom,
  D.~Meyer, T.~Rowley, and S.~M. Robbins.
\newblock Geomview, version 1.8.1, March 2001.
\newblock \burl{http://www.geomview.org}.

\bibitem{mp-rpdgp-93}
W.~Moser and J.~Pach.
\newblock Research problems in discrete geometry: Packing and covering.
\newblock Tech. Rep. 93--32, DIMACS, 1993.

\bibitem{m-ccp-57}
T.~S. Motzkin.
\newblock Comonotone curves and polyhedra.
\newblock \emph{Bull. Amer. Math. Soc.} 63:35, 1957.

\bibitem{r-sphis-47}
R.~Rado.
\newblock A sequence of polyhedra having intersections of specified dimensions.
\newblock \emph{J. London Math. Soc.} 22:287--289, 1947.

\bibitem{s-ravcp-63+}
\textcyr{Yu}.~\textcyr{A}. {\textcyr{Shashkin} [Yu. A. {\v S}a{\v s}kin]}.
\newblock \textcyr{Zamechanie o sosednikh vershniakh na vypuklom
  mnogogranannike} [{A} remark on adjacent vertices on a convex polyhedron].
\newblock \emph{\textcyr{\em Uspekhi Matem. Nauk} [Uspehi Mat. Nauk]}
  18(5):209--211, 1963.

\bibitem{s-tcp-68}
G.~C. Shepherd.
\newblock A theorem on cyclic polytopes.
\newblock \emph{Israel J. Math.} 6: 368--372, 1968.

\bibitem{s-eubnf-91}
R.~Seidel.
\newblock Exact upper bounds for the number of faces in $d$-dimensional
  {V}oronoi diagrams.
\newblock \emph{Applied Geometry and Discrete Mathematics: The Victor Klee
  Festschrift}, pp. 517--530. DIMACS Series in Discrete Mathematics and
  Theoretical Computer Science~4, AMS Press, 1991.

\bibitem{s--97}
P.~St{\"a}ckel.
\newblock \emph{Z. Math. Phys.} 42:275, 1897.

\bibitem{t-udpdn-05}
H.~Tietze.
\newblock {\"U}ber das {P}roblem der {N}achbargite im {R}aum.
\newblock \emph{Monatshefte f{\"u}r Mathematik und Physik} 16:211--216, 1905.

\bibitem{t-fpmsu-65}
H.~Tietze.
\newblock \emph{Famous Problems of Mathematics: Solved and Unsolved
  Mathematical Problems from Antiquity to Modern Times}.
\newblock Graylock Press, New York, 1965.
\newblock Translation of \emph{Gel{\"o}ste und ungel{\"o}ste mathematische
  Probleme aus alter und neuer Zeit}, 1959.

\bibitem{z-alnfs-86}
J.~Zaks.
\newblock Arbitrarily large neighborly families of symmetric convex polytopes.
\newblock \emph{Geom. Dedicata} 20(2):175--179, 1986.

\bibitem{z-nfccp-87}
J.~Zaks.
\newblock Neighborly families of congruent convex polytopes.
\newblock \emph{Amer. Math. Monthly} 94:151--155, 1987.

\bibitem{z-pc-01}
J.~Zaks.
\newblock Personal communication, April 2001.

\bibitem{z-lp-94}
G.~M. Ziegler.
\newblock \emph{Lectures on Polytopes}.
\newblock Graduate Texts in Mathematics 152. Springer-Verlag, Heidelberg, 1994.

\end{thebibliography}

\end{document}